\documentclass[11pt,oneside]{amsart}
\usepackage{mathtools,amssymb,amsthm}
\usepackage{microtype}
\usepackage[hidelinks]{hyperref}
\usepackage[margin=1.05in]{geometry}
\allowdisplaybreaks[2]
\newtheorem{mainthm}{Theorem}

\newtheorem{theorem}{Theorem}[section]
\newtheorem{lemma}[theorem]{Lemma}

\newtheorem{corollary}[theorem]{Corollary}
\theoremstyle{remark}

\numberwithin{equation}{section}
\newcommand{\C}{\mathbb C}
\newcommand{\R}{\mathbb R}
\newcommand{\D}{\mathbb D}
\newcommand{\T}{\mathbb T}
\newcommand{\Rea}{\operatorname{Re}}
\newcommand{\Arg}{\operatorname{Arg}}
\newcommand{\norm}[1]{\lVert #1\rVert}

\title[Sharp weighted Riesz inequalities]{Sharp weighted Riesz inequalities for harmonic mappings in the range $1<p<2$}
\author{David Kalaj}
\address{Faculty of Natural Sciences and Mathematics, University of Montenegro, Cetinjski put b.b., 81000 Podgorica, Montenegro}
\email{davidk@ucg.ac.me}
\author{Anton Gjokaj}
\address{Faculty of Natural Sciences and Mathematics, University of Montenegro, Cetinjski put b.b., 81000 Podgorica, Montenegro}
\email{antondj@ucg.ac.me}

\begin{document}
\begin{abstract}
For $1<p<2$ and $c>0$, we determine the exact norm of the weighted analytic--coanalytic projection operator
$F\mapsto(P_+F,\sqrt c\,P_-F)$ from $L^p(\T)$ to $L^p(\T;\ell_2^2)$.
The sharp constant is
\[
 a_{p,c}=\left(\frac{1+c+\sqrt{1+c^2+2c\cos(2\pi/p)}}{2\sin^2(\pi/p)}\right)^{1/2}.
\]
We construct an explicit plurisubharmonic minorant of the corresponding pointwise obstacle and obtain an origin-corrected inequality for arbitrary harmonic decompositions. The uncorrected inequality holds under an optimal angular condition at the origin. The scalar estimate underlying the construction is proved by angular minimization and a two-sided comparison with solutions of $y''=y$. An explicit family establishes sharpness.
\end{abstract}
\maketitle

\section{Introduction and main results}\label{sec:intro}

Let $\D$ be the unit disk, let $\T=\partial\D$, and write
$dm(e^{it})=dt/(2\pi)$.
For a complex-valued harmonic mapping $f$ on $\D$, set
\[
 \norm{f}_{h^p}
 =\sup_{0<r<1}
 \left(\int_\T |f(r\zeta)|^p\,dm(\zeta)\right)^{1/p}.
\]
We consider decompositions $f=g+\overline h$, where $g$
and $h$ are holomorphic. The boundary $L^p$ norm and
the harmonic Hardy norm agree when $f\in h^p$.

Sharp inequalities relating a harmonic mapping to its
analytic and coanalytic parts have a long history,
originating in the classical conjugate-function
inequalities. Pichorides \cite{Pichorides1972}
determined the best constants in these inequalities,
while Verbitsky \cite{Verbitsky1984} established
sharp Hardy-space estimates in terms of the real and
imaginary parts. Hollenbeck and Verbitsky
\cite{HV2000} determined the exact norm of the
Riesz projection and subsequently investigated
simultaneous estimates for the analytic and
coanalytic projections \cite{HV2010}.

Sharp simultaneous estimates for the analytic and
coanalytic Riesz projections have also been studied
with a general exponent $s>0$.
Melentijevi\'c and Markovi\'c \cite{MM2023}
obtained sharp estimates for $p\ge2$ and
$0<s\le p$. Melentijevi\'c
\cite{Melentijevic2024} subsequently settled
the Hollenbeck--Verbitsky conjecture for
$p\ge2$ and all $s>0$, and established
further sharp results for $1<p\le4/3$ and
$0<s\le\sec^2(\pi/(2p))$.

More recently, Jaguzovi\'c \cite{Jaguzovic2026}
obtained sharp estimates for $4/3<p<2$ and
$0<s\le p/(p-1)$, including the quadratic
case $s=2$. The sharp unweighted quadratic
inequality for harmonic mappings had already
been established by Kalaj \cite{Kalaj2019}
for the entire range $1<p<\infty$.

The corresponding weighted problem is different. 
Here $c>0$ is a constant weight on the
coanalytic component, while the underlying
boundary measure remains unchanged.
For $p\ge2$, Gjokaj, Kalaj, and Vujadinovi\'c
\cite{GKV2025} determined the sharp constant
for the expression
\[
 \left(|g|^2+c|h|^2\right)^{1/2},
 \qquad c>0.
\]
A weighted inequality in the reverse direction
was studied by Kalaj \cite{Kalaj2025}.
The latter result concerns a different operator
estimate and does not give the upper bound
considered here.

The purpose of the present paper is to determine
the exact sharp constant for the weighted quadratic
expression in the subquadratic range $1<p<2$.
More precisely, we allow an arbitrary component
weight $c>0$ and prove that the same explicit
formula as in \cite{GKV2025} remains valid.
Consequently, the weighted analytic--coanalytic
projection norm admits a unified expression
throughout the entire range $1<p<\infty$.
Unlike the previously established unweighted
estimates, the present result captures the
dependence of the optimal constant on the
asymmetric weight $c$.

We also establish a refined inequality for
arbitrary harmonic decompositions
$f=g+\overline h$, involving an explicit
correction term depending on $g(0)h(0)$.
Such a correction is necessary in general.
Indeed, taking the constant functions
$g\equiv1$ and $h\equiv-1$ gives
$f\equiv0$, whereas
$(|g|^2+c|h|^2)^{1/2}>0$.
We identify an optimal angular condition at
the origin under which the correction term
is nonnegative and can therefore be omitted.
The resulting constant remains sharp even under
the stronger normalization $g(0)=h(0)=0$.

Our proof is based on an explicit
plurisubharmonic minorant of the corresponding
pointwise obstacle. The main technical ingredient
is a global hyperbolic support inequality,
established by angular minimization and a
two-sided comparison with solutions of
$y''=y$. The two-sided nature of this comparison
is essential: the logarithmic modulus ratio
$\log(|g|/|h|)$ is unrestricted, and the
contact parameter need not be nonnegative.
An explicit family of harmonic mappings
establishes sharpness of the constant.

Throughout the paper, unless otherwise stated, put
\begin{equation}\label{eq:parameters}
 q=\frac p2,\qquad \tau=\frac{\pi}{p},\qquad
 \kappa=-\cos\tau\in(0,1),\qquad \sigma=\sin\tau>0.
\end{equation}
For $c>0$, define
\begin{equation}\label{eq:sharp-constant}
 \Delta=\sqrt{(1-c)^2+4c\kappa^2},\qquad
 a_{p,c}^2=\frac{1+c+\Delta}{2\sigma^2},\qquad
 R=\frac{1-c+\Delta}{2\kappa},
\end{equation}
and
\begin{equation}\label{eq:minorant-coefficient}
 b_{p,c}=\frac{1+c+\Delta}{\sigma}
           \left(\frac{\Delta}{\kappa}\right)^{q-1}>0.
\end{equation}
The angular function needed in this range is particularly simple:
\begin{equation}\label{eq:G-definition}
 G_p(z,w)=
 \begin{cases}
 |zw|^q\cos\bigl(q\Arg(zw)\bigr),&zw\ne0,\\
 0,&zw=0.
 \end{cases}
\end{equation}
Here $\Arg$ is the principal argument in $[-\pi,\pi]$; the two values at the negative real axis give the same cosine. It is important to take the principal argument of the product, rather than an unwrapped sum of two independently chosen arguments.

\begin{mainthm}[Sharp weighted Riesz inequality]\label{thm:main}
Let $1<p<2$, $c>0$, and let $f=g+\overline h\in h^p$. Then $g,h\in H^p$, and
\begin{equation}\label{eq:corrected-integral}
 \int_\T (|g|^2+c|h|^2)^q\,dm
 +b_{p,c}G_p(g(0),h(0))
 \le a_{p,c}^p\int_\T |g+\overline h|^p\,dm.
\end{equation}
Consequently, if
\begin{equation}\label{eq:sector-condition}
 g(0)h(0)=0\quad\text{or}\quad
 |\Arg(g(0)h(0))|\le\frac{\pi}{p},
\end{equation}
then
\begin{equation}\label{eq:main-norm}
 \left(\int_\T (|g|^2+c|h|^2)^{p/2}\,dm\right)^{1/p}
 \le a_{p,c}\norm{f}_{h^p}.
\end{equation}
The constant is optimal even under the stronger normalization $g(0)=h(0)=0$.
The angular condition is optimal in the following precise sense: for every $\theta\in[-\pi,\pi]$ with $|\theta|>\pi/p$, there are constant functions $g,h$ with $\Arg(g(0)h(0))=\theta$ for which \eqref{eq:main-norm} fails.
\end{mainthm}

The hypothesis $\Rea(g(0)h(0))=0$ used in the companion result is a special case of \eqref{eq:sector-condition}; so is $\Rea(g(0)h(0))\ge0$.
For $F(e^{it})=\sum_{n\in\mathbb Z}\widehat F(n)e^{int}$, use the convention
\[
 P_+F=\sum_{n\ge0}\widehat F(n)e^{int},\qquad
 P_-F=\sum_{n<0}\widehat F(n)e^{int}.
\]
Thus the constant Fourier coefficient belongs to $P_+$.

\begin{corollary}[Full-range weighted projection norm]
\label{cor:projection}
For every $1<p<\infty$ and $c>0$,
\begin{equation}\label{eq:full-range}
\begin{aligned}
&\left\|F\longmapsto
(P_+F,\sqrt c\,P_-F)\right\|
_{L^p(\T)\to L^p(\T;\ell_2^2)} =
\left(
\frac{1+c+\sqrt{1+c^2+2c\cos(2\pi/p)}}
{2\sin^2(\pi/p)}
\right)^{1/2}.
\end{aligned}
\end{equation}
\end{corollary}

For $1<p<2$, the assertion follows from
Theorem~\ref{thm:main}, as proved below.
For $2\le p<\infty$, it follows from the
sharp weighted inequality established in
\cite{GKV2025}. Thus the same explicit
formula determines the operator norm
throughout the entire range $1<p<\infty$.

The mechanism behind Theorem~\ref{thm:main} is the following pointwise result.

\begin{mainthm}[Explicit plurisubharmonic minorant]\label{thm:pointwise}
The function $G_p$ in \eqref{eq:G-definition} is continuous and plurisubharmonic on $\C^2$. For every $z,w\in\C$,
\begin{equation}\label{eq:pointwise}
 (|z|^2+c|w|^2)^q+b_{p,c}G_p(z,w)
 \le a_{p,c}^p|z+\overline w|^p.
\end{equation}
Apart from the origin, equality holds exactly for those $z,w\in\C$
with $zw\ne0$ satisfying
\begin{equation}\label{eq:contact}
 |z|=R|w|,\qquad |\Arg(zw)|=\frac{\pi}{p}.
\end{equation}

For the fixed coefficient $a_{p,c}^p$ and the fixed function $G_p$, the coefficient $b_{p,c}$ is the unique real number for which \eqref{eq:pointwise} can hold globally.
\end{mainthm}

We emphasize that Theorem~\ref{thm:pointwise}
provides an explicit plurisubharmonic minorant
of the corresponding obstacle, without claiming
that it coincides with its plurisubharmonic
envelope.

The paper is organized as follows.
Section~\ref{sec:scalar} establishes the global
hyperbolic support inequality by angular
minimization and a two-sided differential comparison.
Section~\ref{sec:minorant} uses this scalar estimate
to construct the plurisubharmonic minorant and
derive the integral inequality.
Section~\ref{sec:sharpness} proves sharpness
by an explicit extremizing family and discusses
the optimality of the angular condition and
the endpoint cases.

\section{A global hyperbolic support inequality}\label{sec:scalar}
The following strengthening of the one-sided main inequality is the scalar result used in the proof.

\begin{lemma}[Global hyperbolic support]\label{lem:global-support}
Let $1<p<2$, let $u,\xi\in\R$, and suppose that
\[
 V(u,\xi):=\cosh u-\kappa\cosh(u-\xi)\ge0.
\]
Then, for every $t\in[0,\pi]$,
\begin{equation}\label{eq:global-support}
 (\cosh u+\cos t)^q-V(u,\xi)^q
 \ge \sigma(\cosh\xi-\kappa)^{q-1}\cos(qt).
\end{equation}
Equality holds if and only if $u=\xi$ and $t=\tau$.
\end{lemma}

We first prove the trigonometric positivity needed for the differential comparison.

\begin{lemma}[Auxiliary positivity]\label{lem:auxiliary}
Let $1/2<q<1$, $s=1-q$, and $0<t<\pi$. Put
\[
 x=\cos t,\qquad
 v=\frac{\sin(st)}{\sin(qt)},\qquad
 C=v^2-s^2(1+2xv+v^2),
\]
and
\[
 B=2C+v(vx-qv^2-s).
\]
Then $C\ge0$ and $B>0$.
\end{lemma}

\begin{proof}
Since $0<s<q<1$ and the function
$z\mapsto \sin z/z$ is strictly decreasing
on $(0,\pi)$, we have
\[
 v=\frac{\sin(st)}{\sin(qt)}>\frac{s}{q}.
\]
Moreover,
\[
 \sin(qt)-\sin(st)
 =2\cos(t/2)\sin((q-s)t/2)>0,
\]
and consequently
\begin{equation}\label{eq:v-bounds}
 \frac{s}{q}<v<1.
\end{equation}
Using the identity
\[
 1+2xv+v^2
 =\frac{\sin^2t}{\sin^2(qt)}
\]
and the concavity inequality
$\sin(st)\ge s\sin t$, we obtain
\[
 C=
 \frac{\sin^2(st)-s^2\sin^2t}
 {\sin^2(qt)}
 \ge0.
\]
It remains to prove that $B>0$.
We first observe that, as $t\to0$,
\[
 v(t)=\frac{s}{q}
 \left(1+\frac{q^2-s^2}{6}t^2
 +O(t^4)\right),
 \qquad
 x(t)=1-\frac{t^2}{2}+O(t^4).
\]
Substituting these expansions into the
definition of $B$ gives
\begin{equation}\label{eq:B-small}
 B(t)=\frac{s^2(q-s)}{3q^2}t^2
 +O(t^4).
\end{equation}
Since $q>s$, it follows that
$B(t)>0$ for all sufficiently small $t>0$.
Differentiating
$v=\sin(st)/\sin(qt)$ and $x=\cos t$,
we obtain
\[
 \sin t\,v'(t)
 =N:=s+(s-q)xv-qv^2,
\]
and
\[
 \sin t\,x'(t)=-(1-x^2).
\]
An expansion of the definition of $B$ gives
\begin{equation}\label{eq:B-expanded}
 B=v(v-4s^2)x
   +2(1-s^2)v^2-qv^3-sv-2s^2.
\end{equation}
Notice that the coefficient of $x$
is strictly positive. Indeed,
by \eqref{eq:v-bounds},
\[
 v-4s^2>
 \frac{s}{q}-4s^2
 =\frac{s(1-4sq)}{q}
 =\frac{s(q-s)^2}{q}>0.
\]
We now show that $B$ cannot vanish
in $(0,\pi)$. Suppose that $B(t)=0$
at some $t\in(0,\pi)$.
Solving \eqref{eq:B-expanded} for
$x=\cos t$, we find
\begin{equation}\label{eq:x-at-zero}
 x=
 \frac{qv^3+sv+2s^2-2(1-s^2)v^2}
      {v(v-4s^2)}.
\end{equation}
The denominator is strictly positive
by the preceding observation.
Regarding $B$ as a polynomial
$B=B(x,v)$, its partial derivatives are
\[
 B_x=v(v-4s^2),
\]
and
\[
 B_v=(2v-4s^2)x
      +4(1-s^2)v-3qv^2-s.
\]
Hence, by the chain rule,
\begin{equation}\label{eq:B-chain-rule}
 \sin t\,B'(t)
 =B_vN-B_x(1-x^2).
\end{equation}
At a zero of $B$, substitute
\eqref{eq:x-at-zero} into
\eqref{eq:B-chain-rule}.
After clearing denominators and
collecting terms, we obtain
\begin{equation}\label{eq:B-zero-derivative}
 \sin t\,B'(t)
 =
 \frac{(1-v)(qv+s)(qv-s)}
 {(v-4s^2)^2}\,Q_s(v),
\end{equation}
where, writing $\rho=qv-s$,
\begin{equation}\label{eq:Q-positive}
\begin{aligned}
(q-s)Q_s(v)
={}&s^2(5-6s)(1-v)^3\\
&+4s(1-2s+s^2+2s^3)
       \rho(1-v)^2\\
&+(1+2s)(3-8s^2)
       \rho^2(1-v)\\
&+4s(1+2s)\rho^3.
\end{aligned}
\end{equation}

We claim that $Q_s(v)>0$.
Indeed, since $0<s<1/2$, we have
\[
 5-6s>2,\qquad
 3-8s^2>1,
\]
and
\[
 1-2s+s^2+2s^3
 =(1-s)^2+2s^3>0.
\]
Furthermore, \eqref{eq:v-bounds} yields
\[
 \rho=qv-s>0,\qquad 1-v>0.
\]
Thus every term on the right-hand
side of \eqref{eq:Q-positive} is strictly
positive. Since $q-s>0$, this proves
\[
 Q_s(v)>0.
\]
All the remaining factors on the
right-hand side of
\eqref{eq:B-zero-derivative} are
strictly positive. Therefore,
\[
 B'(t)>0
 \quad\text{whenever }B(t)=0,
 \qquad 0<t<\pi.
\]
Finally, suppose that $B$ has a zero
in $(0,\pi)$, and let $t_0$ be its
first zero. Such a first zero exists
because $B$ is continuous and,
by \eqref{eq:B-small}, strictly
positive near the origin.
Then $B(t)>0$ for $0<t<t_0$,
which implies $B'(t_0)\le0$.
This contradicts the preceding
conclusion that $B'(t_0)>0$.
Consequently,
\[
 B(t)>0,\qquad 0<t<\pi,
\]
as required.
\end{proof}

\subsection{Angular minimization and its differential inequality}
\begin{proof}[Proof of Lemma~\ref{lem:global-support}]
Write $s=1-q$. We first fix $b\ge0$ and define, for \emph{every} $u\ge0$,
\begin{equation}\label{eq:F-m-definition}
 \begin{gathered}
 D=\cosh b-\kappa,\qquad K=\sigma D^{-s},\\
 F(u,t)=(\cosh u+\cos t)^q-K\cos(qt),\qquad
 m(u)=\min_{0\le t\le\pi}F(u,t).
 \end{gathered}
\end{equation}
The range $0\le u<b$ is essential, since in the final reduction
$a=|u|$, $b=|\xi|$, one does not necessarily have $a\ge b$.
\\
For $A\ge1$ and $0<t<\pi$, put
\[
 R_A(t)=\frac{\sin t}{\sin(qt)(A+\cos t)^s}.
\]
Here $R_A$ is an auxiliary angular function, unrelated to the contact radius $R$ in \eqref{eq:sharp-constant}.
With $x=\cos t$ and $v=\sin(st)/\sin(qt)$, direct differentiation gives
\begin{equation}\label{eq:R-log-derivative}
 \frac{d}{dt}\log R_A(t)=-\frac{J}{\sin t\,(A+\cos t)},
\end{equation}
where
\begin{equation}\label{eq:J-definition}
\begin{aligned}
 J&=(qv-sx)A+qvx-s\\
  &=(A-1)(qv-sx)+(1+x)(qv-s)>0.
\end{aligned}
\end{equation}
Indeed, $qv-s>0$ and $qv-sx=(qv-s)+s(1-x)>0$.
Thus $R_A$ is strictly decreasing. Its continuous endpoint values are
\[
 R_A(0)=\frac{(A+1)^{-s}}q,\qquad R_A(\pi)=0.
\]
For $A=1$, the latter limit uses $1-2s=2q-1>0$.
Since
\begin{equation}\label{eq:F-angle-derivative}
 F_t(u,t)=q\sin(qt)\bigl(K-R_{\cosh u}(t)\bigr),
\end{equation}
there is a unique minimizing angle $\theta(u)$. If
$K<R_{\cosh u}(0)$, it is the unique interior solution of
\begin{equation}\label{eq:stationarity}
 R_{\cosh u}(\theta(u))=K;
\end{equation}
otherwise the minimum is attained at $t=0$, so $\theta(u)=0$.
At $u=b$,
\[
 R_{\cosh b}(\tau)=K,
\]
so $\theta(b)=\tau$. Since $\tau>0$ and $R_{\cosh b}$ is strictly
decreasing, $R_{\cosh b}(0)>K$. As
$u\mapsto R_{\cosh u}(0)$ decreases strictly to zero, there is exactly
one transition value $u_*>b$ such that $\theta(u)>0$ for
$0\le u<u_*$ and $\theta(u)=0$ for $u\ge u_*$.
\\
We need positivity of $m$ also to the left of $b$.
At any interior minimum, the identity
\[
 \sin\theta\cot(q\theta)=\cos\theta+
                \frac{\sin(s\theta)}{\sin(q\theta)}
\]
and \eqref{eq:stationarity} give
\begin{equation}\label{eq:m-positive}
 m(u)=(A+x)^{-s}(A-v)>0,
 \quad A=\cosh u,\quad x=\cos\theta(u),\quad
 v=\frac{\sin(s\theta(u))}{\sin(q\theta(u))}.
\end{equation}
The strict inequality follows from $A\ge1>v$. In particular $m(0)>0$.
For fixed $t$, $F(u,t)$ is nondecreasing in $u$, so
$m(u)\ge m(0)>0$ for all $u\ge0$.
We may therefore define
\[
 h(u)=m(u)^{1/q},\qquad u\ge0.
\]
At $u=b$, using $\theta(b)=\tau$, $q\tau=\pi/2$, and
$\sin(s\tau)=\kappa$, differentiation at the minimum yields
\begin{equation}\label{eq:initial-data}
 h(b)=D,\qquad h'(b)=\sinh b.
\end{equation}
For $b=0$, the derivative is the right derivative.
\\
On the interior branch, abbreviate $S=\sinh u$ and retain the notation
of \eqref{eq:m-positive}. Since $m(u)=F(u,\theta(u))$ and
$F_t(u,\theta(u))=0$, differentiation at the minimum gives
\begin{equation}\label{eq:h-log-derivative}
 m'(u)=qS(A+x)^{-s},\qquad \frac{h'}h=\frac{S}{A-v}.
\end{equation}
The identity for $v'$ in the proof of Lemma~\ref{lem:auxiliary}, and
implicit differentiation in \eqref{eq:stationarity}, give
\[
 \theta'(u)=-\frac{sS\sin\theta}{J},\qquad
 \frac{dv}{du}=-\frac{sSN}{J},\qquad
 N=s+(s-q)xv-qv^2.
\]
It follows that
\begin{equation}\label{eq:h-ode}
 \frac{h''-h}{h}
 =\frac{v(A-v)J-s(A^2-1)N}{(A-v)^2J}.
\end{equation}
For the quantities $C,B$ in Lemma~\ref{lem:auxiliary}, the numerator
has the exact decomposition
\begin{equation}\label{eq:numerator-decomposition}
\begin{aligned}
 v(A-v)J-s(A^2-1)N
 ={}&C(A-1)^2+B(A-1)\\
    &+v(1-v)(1+x)(qv-s).
\end{aligned}
\end{equation}
All terms are nonnegative, and the last is strictly positive. Hence
\begin{equation}\label{eq:strict-ode}
 h''>h\qquad(0\le u<u_*),
\end{equation}
with one-sided interpretation at zero. Notice that
Lemma~\ref{lem:auxiliary} applies to every $0<\theta<\pi$, including
the angles $\theta>\tau$ occurring when $u<b$.
\\
For $u>u_*$ the minimum is at zero, and
\[
 m(u)=(A+1)^q-K.
\]
Set $w=K(A+1)^s-1$. The minimum condition gives
$w\ge s/q>0$. Moreover, $m>0$ implies
$K<(A+1)^q$, and hence $w<A$. Thus $A-w>0$, while
\[
 w'=sK(A+1)^{s-1}S\ge0.
\]
Consequently
\begin{equation}\label{eq:boundary-ode}
 \frac{h'}h=\frac{S}{A-w},\qquad
 \frac{h''-h}{h}=\frac{w}{A-w}+\frac{Sw'}{(A-w)^2}>0.
\end{equation}

As $u\uparrow u_*$, the stationarity condition
\eqref{eq:stationarity} yields
$\theta(u)\to0$, and hence $v\to s/q$.
On the other hand, the defining relation
$K=R_{\cosh u_*}(0)$ gives $w(u_*)=s/q$.
The values of $h$ and $h'$ match, so
$h\in C^1([0,\infty))$ and is smooth away from the transition.
To make the distributional assertion precise, the function
\[
 u\longmapsto h'(u)-\int_0^u h(r)\,dr
\]
is increasing on each smooth branch and continuous at $u_*$. It is
therefore increasing globally. Thus
\[
 \mu=h''-h
\]
is a nonnegative Radon measure, with no atom at the transition.
Its density is strictly positive on each open smooth branch.

\medskip
\noindent\emph{Two-sided hyperbolic comparison.}
The solution with the initial data \eqref{eq:initial-data} is
\begin{equation}\label{eq:support-solution}
\begin{aligned}
 y_b(u)&=D\cosh(u-b)+\sinh b\,\sinh(u-b)\\
       &=\cosh u-\kappa\cosh(u-b).
\end{aligned}
\end{equation}
The variation-of-constants formula gives
\begin{equation}\label{eq:two-sided-comparison}
 h(u)-y_b(u)=
 \begin{cases}
 \displaystyle\int_b^u\sinh(u-r)\,d\mu(r),&u\ge b,\\[6pt]
 \displaystyle\int_u^b\sinh(r-u)\,d\mu(r),&0\le u\le b.
 \end{cases}
\end{equation}
Both integrals are nonnegative. If $u\ne b$, the appropriate interval
is nontrivial and contains points of a smooth branch where the density
of $\mu$ is strictly positive; hence the integral is strictly positive.
Therefore $h(u)\ge y_b(u)$ for every $u\ge0$, on both sides of $b$.
Whenever $y_b(u)\ge0$, raising to the power $q$ yields
\[
 F(u,t)\ge m(u)=h(u)^q\ge y_b(u)^q.
\]
Equality requires $u=b$ and then, by uniqueness of the angular
minimum, $t=\tau$.
\\
Finally allow $u,\xi\in\R$ as in the lemma. Put $a=|u|$ and
$b=|\xi|$. We have $|u-\xi|\ge|a-b|$, and consequently
\[
 0\le V(u,\xi)\le\cosh a-\kappa\cosh(a-b)=y_b(a).
\]
The estimate already proved at $(a,b)$ implies
\eqref{eq:global-support}. For equality, the strict comparison forces
$a=b$; equality of the preceding hyperbolic expressions then forces
$u=\xi$. The remaining condition is $t=\tau$.
\end{proof}

\begin{corollary}[The one-sided main inequality]\label{lem:main-inequality}
Let $1<p\le2$, $\alpha\ge0$, and $0\le\gamma\le\alpha$. For every $t\in[0,\pi]$,
\begin{equation}\label{eq:original-main-inequality}
\begin{aligned}
 &(\cosh\alpha+\cos t)^{p/2}
 -(\cosh\alpha+\cosh\gamma\cos(\pi/p))^{p/2}\\
 &\qquad\ge
 \frac{\cos(pt/2)\sin(\pi/p)}
 {(\cosh(\alpha-\gamma)+\cos(\pi/p))^{1-p/2}}.
\end{aligned}
\end{equation}
For $1<p<2$, equality holds exactly when $\gamma=0$ and $t=\pi/p$.
For $p=2$, the formula is an identity.
\end{corollary}
\begin{proof}
For $p<2$, take $u=\alpha$ and $\xi=\alpha-\gamma$ in Lemma~\ref{lem:global-support}.
The base is positive because
$\cosh\alpha-\kappa\cosh\gamma\ge(1-\kappa)\cosh\alpha>0$.
The endpoint $p=2$ is immediate.
\end{proof}

\section{The minorant and the integral inequality}\label{sec:minorant}
\subsection{Plurisubharmonicity}
The angular profile is already present in the unweighted construction of \cite{Kalaj2019}. We verify its branch convention here; the weight dependence is carried by \eqref{eq:minorant-coefficient}.
\begin{lemma}[The product-power minorant]\label{lem:psh}
For $1/2<q<1$, the continuous function
\[
 \Psi_q(\zeta)=|\zeta|^q\cos(q\Arg\zeta),\qquad \Psi_q(0)=0,
\]
is subharmonic on $\C$. Consequently $G_p(z,w)=\Psi_q(zw)$ is plurisubharmonic on $\C^2$.
\end{lemma}
\begin{proof}
Away from the negative real axis and the origin, $\Psi_q$ is locally
the real part of a holomorphic branch of $\zeta^q$, and is therefore
harmonic.
Consider the $2\pi$-periodic extension of the angular profile
\[
 \phi(\theta)=\cos(q\theta),\qquad -\pi\le\theta\le\pi.
\]
The endpoint values agree, while the derivative has the jump
\[
 \phi'(-\pi+)-\phi'(\pi-)=2q\sin(q\pi)>0.
\]
Since $\phi''+q^2\phi=0$ on $(-\pi,\pi)$, it follows in the periodic
distributional sense that
\[
 \phi''+q^2\phi
 =2q\sin(q\pi)\,\delta_{\pi},
\]
which is a positive measure. More explicitly, in the plane,
\[
 \langle\Delta\Psi_q,\chi\rangle
 =2q\sin(q\pi)\int_0^\infty r^{q-1}\chi(-r)\,dr,
 \qquad \chi\in C_c^\infty(\C).
\]
The integral is locally finite at zero. There is no additional mass
there: the boundary terms from a circle of radius $\varepsilon$ are
$O(\varepsilon^q)$ and tend to zero. Thus $\Delta\Psi_q$ is a
nonnegative measure, and $\Psi_q$ is subharmonic on $\C$.\\
Finally, $(z,w)\mapsto zw$ is holomorphic. Its restriction to every
complex line is holomorphic, and the composition of a subharmonic
function with a holomorphic map is subharmonic. Hence
$G_p(z,w)=\Psi_q(zw)$ is subharmonic on every complex line, and
therefore plurisubharmonic on $\C^2$.
\end{proof}

\subsection{The weight and its contact radius}
\begin{lemma}[Contact identities]\label{lem:contact-identities}
The constants in \eqref{eq:sharp-constant} satisfy
\begin{equation}\label{eq:contact-identities}
 \kappa<R<\kappa^{-1},\qquad
 c=\frac{R(1-\kappa R)}{R-\kappa},\qquad
 a_{p,c}^2=\frac{R}{R-\kappa}>\max\{1,c\}.
\end{equation}
If $\xi=\log R$ and $D_R=\cosh\xi-\kappa$, then
\begin{equation}\label{eq:b-contact-form}
 b_{p,c}=2^q a_{p,c}^p\sigma D_R^{q-1}.
\end{equation}
For every $r>0$, writing $u=\log r$,
\begin{equation}\label{eq:weight-reduction}
 \frac{r^2+c}{2r a_{p,c}^2}
 =\cosh u-\kappa\cosh(u-\xi)>0.
\end{equation}
\end{lemma}
\begin{proof}
The positive root of
\[
 \kappa R^2-(1-c)R-c\kappa=0
\]
is the number in \eqref{eq:sharp-constant}. Since
\[
 P(\kappa)=\kappa(\kappa^2-1)<0,\qquad
 P(\kappa^{-1})=c(\kappa^{-1}-\kappa)>0,
\]
we get $\kappa<R<\kappa^{-1}$. Solving for $c$ gives the second identity in \eqref{eq:contact-identities}.
The number $\lambda=R/(R-\kappa)$ satisfies
\[
 (1-\kappa^2)\lambda^2-(1+c)\lambda+c=0,
\]
and
\[
 \lambda-1=\frac{\kappa}{R-\kappa}>0,\qquad
 \lambda-c=\frac{\kappa R^2}{R-\kappa}>0.
\]
Hence $\lambda$ is the larger root, so $\lambda=a_{p,c}^2$.
If $\xi=\log R$, then
\[
 D_R=\cosh\xi-\kappa=\frac{R^2+1-2\kappa R}{2R}
 =\frac{\Delta}{2\kappa a_{p,c}^2},
\]
where the last identity follows from
\[
 \Delta=\frac{\kappa(R^2+1-2\kappa R)}{R-\kappa}.
\]
Substitution gives \eqref{eq:b-contact-form}.
Finally, for $u=\log r$,
\[
 \cosh u-\kappa\cosh(u-\xi)
 =\frac12\left[r\left(1-\frac{\kappa}{R}\right)
 +\frac1r(1-\kappa R)\right]
 =\frac{r^2+c}{2r a_{p,c}^2}>0,
\]
by \eqref{eq:contact-identities}.
\end{proof}

\begin{proof}[Proof of Theorem~\ref{thm:pointwise}]
Plurisubharmonicity follows from Lemma~\ref{lem:psh}.
Assume $zw\ne0$, and put
\[
 r=\frac{|z|}{|w|},\qquad u=\log r,\qquad
 \xi=\log R,\qquad t=|\Arg(zw)|.
\]
Then
\[
 |z+\overline w|^2=2|zw|(\cosh u+\cos t),\qquad
 |z|^2+c|w|^2=2|zw|\frac{r^2+c}{2r}.
\]
Apply Lemma~\ref{lem:global-support}, using \eqref{eq:weight-reduction}, and multiply by $a_{p,c}^p(2|zw|)^q$.
Formula \eqref{eq:b-contact-form} gives exactly \eqref{eq:pointwise}.
The equality conditions are $u=\xi$ and $t=\tau$, which are \eqref{eq:contact}.
If precisely one of $z,w$ vanishes, the inequality is strict by $a_{p,c}^2>\max\{1,c\}$; the origin is an equality point.

To verify the uniqueness of the correction coefficient, suppose that \eqref{eq:pointwise} holds with a real number $b$ in place of $b_{p,c}$.
For $t\in[0,\pi]$, take $z=R e^{it/2}$, $w=e^{it/2}$ and define
\[
 E(t)=a_{p,c}^p(R^2+1+2R\cos t)^q-(R^2+c)^q-bR^q\cos(qt).
\]
The contact identities, together with $q\tau=\pi/2$, give $E(\tau)=0$, independently of $b$.
Since $E(t)\ge0$ on $[0,\pi]$ and $\tau\in(0,\pi)$, the point $t=\tau$ is an interior minimum of $E$, hence $E'(\tau)=0$.
Therefore
\[
 b=2a_{p,c}^p\sigma R^{1-q}(R^2+1-2\kappa R)^{q-1}
   =b_{p,c},
\]
as claimed.
\end{proof}

\subsection{Integration and the condition at the origin}
\begin{proof}[Upper bounds in Theorem~\ref{thm:main} and Corollary~\ref{cor:projection}]
We establish the upper bounds in the subquadratic
range $1<p<2$, which is the only case requiring
proof in Corollary~\ref{cor:projection}.
For fixed $0<r<1$, apply \eqref{eq:pointwise} to $z=g(r\zeta)$ and $w=h(r\zeta)$.
Since $G_p(g,h)$ is subharmonic,
\[
 \int_\T G_p(g(r\zeta),h(r\zeta))\,dm(\zeta)
 \ge G_p(g(0),h(0)).
\]
After integration, this yields
\begin{equation}\label{eq:radial-corrected}
\begin{aligned}
 &\int_\T (|g(r\zeta)|^2+c|h(r\zeta)|^2)^q\,dm(\zeta)
 +b_{p,c}G_p(g(0),h(0))\\
 &\qquad\le a_{p,c}^p\int_\T |f(r\zeta)|^p\,dm(\zeta).
\end{aligned}
\end{equation}
The right side is bounded uniformly in $r$. Since the correction is a fixed finite number, \eqref{eq:radial-corrected} bounds the $H^p$ means of both $g$ and $h$.
Thus $g,h\in H^p$, even without a sign assumption on the correction.
Their radial functions converge in $L^p$ to their boundary values.
Letting $r\uparrow1$ gives \eqref{eq:corrected-integral}.
By \eqref{eq:G-definition}, condition \eqref{eq:sector-condition} is exactly
$G_p(g(0),h(0))\ge0$, and dropping the correction gives \eqref{eq:main-norm}.

For the projection formulation, take the Poisson extension of $F\in L^p(\T)$ and its decomposition normalized by $h(0)=0$.
Its boundary analytic and coanalytic parts are $P_+F$ and $P_-F=\overline h$.
The upper bound in Corollary~\ref{cor:projection} follows.
Sharpness and optimality of the angular condition are proved below.
\end{proof}

\section{Sharpness and consequences}\label{sec:sharpness}

\subsection{An extremizing family}
\begin{proof}[Sharpness in Theorem~\ref{thm:main} and Corollary~\ref{cor:projection}]
For $0<\beta<1/p$, use the holomorphic branch on $\D$ to define
\begin{equation}\label{eq:extremizing-family}
 W_\beta(z)=(1+z)^{-\beta}-1,\qquad
 g_\beta=R W_\beta,\qquad h_\beta=W_\beta,\qquad
 f_\beta=g_\beta+\overline{h_\beta}.
\end{equation}
Since $p\beta<1$, the boundary singularity $|1+e^{it}|^{-p\beta}$ is integrable near $-1$, hence $W_\beta\in H^p$, and $g_\beta(0)=h_\beta(0)=0$.
Let
\[
 N_\beta=\int_\T |W_\beta|^p\,dm.
\]
As $\beta\uparrow1/p$, $N_\beta\to\infty$.
Indeed, near $-1$ the leading singularity is $|1+e^{it}|^{-p\beta}$; subtracting $1$ does not alter its divergent integral in the limit.
On every closed arc avoiding $-1$, $W_\beta$ remains uniformly bounded.
Define the probability measures
\[
 d\nu_\beta
 =\frac{|W_\beta|^p}{N_\beta}\,dm.
\]
For every open neighborhood $U$ of $-1$
in $\T$, the functions $W_\beta$ are
uniformly bounded on $\T\setminus U$
as $\beta\uparrow1/p$. Hence
\[
 \nu_\beta(\T\setminus U)
 =\frac{1}{N_\beta}
  \int_{\T\setminus U}|W_\beta|^p
  \,dm
 \le\frac{C_U}{N_\beta}
 \longrightarrow0.
\]
Thus the measures $\nu_\beta$
concentrate at $-1$ as
$\beta\uparrow1/p$.
\\
For $-\pi<t<\pi$,
\[
 1+e^{it}=2\cos(t/2)e^{it/2}.
\]
Moreover,
\[
 \frac{W_\beta(e^{it})}{(1+e^{it})^{-\beta}}
 =1-(1+e^{it})^\beta\longrightarrow1
\]
near either endpoint, uniformly for $\beta$ in a fixed left neighborhood of $1/p$.
Consequently, as $(\beta,t)\to(1/p,\pi-)$ or $(1/p,-\pi+)$, respectively,
\[
 \frac{W_\beta(e^{it})}{|W_\beta(e^{it})|}
 \longrightarrow e^{-i\tau/2}\quad\text{or}\quad e^{i\tau/2}.
\]
Hence the corresponding values of
$|R W_\beta+\overline{W_\beta}|^p/|W_\beta|^p$
have the same limit,
\[
 \left|Re^{-i\tau/2}+e^{i\tau/2}\right|^p
 =(R^2+1+2R\cos\tau)^q,
\]
and this quotient is bounded by $(R+1)^p$.
Concentration therefore gives
\begin{equation}\label{eq:sharpness-limit}
 \lim_{\beta\uparrow1/p}
 \frac{\int_\T|f_\beta|^p\,dm}{N_\beta}
 =(R^2+1+2R\cos\tau)^q.
\end{equation}
On the other hand,
\[
 \int_\T(|g_\beta|^2+c|h_\beta|^2)^q\,dm
 =(R^2+c)^qN_\beta.
\]
Using $q=p/2$, $\cos\tau=-\kappa$, and the contact identity
\[
 R^2+c=a_{p,c}^2(R^2+1-2\kappa R),
\]
it follows that
\[
 \lim_{\beta\uparrow1/p}
 \frac{\norm{(|g_\beta|^2+c|h_\beta|^2)^{1/2}}_{L^p}}
 {\norm{f_\beta}_{h^p}}
 =\left(\frac{R^2+c}{R^2+1-2\kappa R}\right)^{1/2}
 =a_{p,c}.
\]
Since $g_\beta(0)=h_\beta(0)=0$ for every $\beta$, this proves sharpness even under the stronger normalization at the origin. For $F_\beta=f_\beta|_\T$, one has $P_+F_\beta=g_\beta$ and $P_-F_\beta=\overline{h_\beta}$, so the same family proves sharpness in Corollary~\ref{cor:projection}.
\end{proof}

The limiting profile is
\[
 W_*(z)=(1+z)^{-1/p}-1,\qquad f_*(z)=R W_*(z)+\overline{W_*(z)}.
\]
The convergence is locally uniform, also for the analytic and coanalytic derivatives.
The nonzero leading boundary coefficient gives
$|f_*(e^{it})|^p\asymp(\pi-|t|)^{-1}$ near $t=\pm\pi$, so $f_*\notin h^p$.
For $c=1$, $R=1$ and the profile is real harmonic.
For $c\ne1$, it is an invertible real-linear image of a conformal map; it is sense-preserving when $R>1$ and sense-reversing when $R<1$.

\subsection{Optimality of the angular normalization}
\begin{proof}[Final assertion of Theorem~\ref{thm:main}]
Fix $\theta\in[-\pi,\pi]$ with $|\theta|>\tau$ and take the constant functions
\[
 g(z)=R e^{i\theta/2},\qquad h(z)=e^{i\theta/2}.
\]
Their product has principal argument $\theta$ (with the usual equivalent choices at $\pm\pi$).
Since
\[
 \cos\theta=\cos|\theta|<\cos\tau=-\kappa,
\]
we have
\[
 R^2+1+2R\cos\theta<R^2+1-2\kappa R.
\]
Using the contact identity,
\[
 R^2+c=a_{p,c}^2(R^2+1-2\kappa R)
       >a_{p,c}^2|g+\overline h|^2.
\]
Thus \eqref{eq:main-norm} fails. This includes the case $g+\overline h=0$.
\end{proof}

\subsection{The unweighted case, the endpoint, and the full-range formula}
For $c=1$, formulas \eqref{eq:sharp-constant}--\eqref{eq:minorant-coefficient} reduce to
\begin{equation}\label{eq:unweighted}
 R=1,\qquad
 a_{p,1}=\frac1{\sqrt{1+\cos(\pi/p)}}
        =\frac1{\sqrt2\cos(\pi/(2p))},\qquad
 b_{p,1}=2^{p/2}\tan\frac{\pi}{2p}.
\end{equation}
These are the unweighted subquadratic constants from \cite{Kalaj2019}.
The formula also gives the symmetry and limit
\[
 a_{p,c}=\sqrt c\,a_{p,1/c},\qquad
 \lim_{c\downarrow0}a_{p,c}=\csc\frac{\pi}{p}.
\]
The latter recovers the upper Riesz-projection bound. Since $R\to1/\kappa$ as $c\downarrow0$, its sharpness can also be read directly from \eqref{eq:extremizing-family}, taking the fixed radius $R=1/\kappa$ and retaining only $g_\beta$ in the numerator. The limit of the ratio is $1/\sigma$, in agreement with \cite{HV2000}.
\\
At $p=2$, the formula extends continuously to
\[
 a_{2,c}^2=\max\{1,c\}.
\]
For every harmonic decomposition $f=g+\overline h\in h^2$,
\[
 \int_\T|f|^2\,dm
 =\int_\T(|g|^2+|h|^2)\,dm+2\Rea(g(0)h(0)).
\]
Moreover,
\[
 |g|^2+c|h|^2\le a_{2,c}^2(|g|^2+|h|^2).
\]
Therefore the sharp uncorrected estimate holds whenever
$\Rea(g(0)h(0))\ge0$, and in particular under the normalization used in the companion paper.
Under $\Rea(g(0)h(0))=0$, equality holds exactly when $h\equiv0$ if $c<1$, when $g\equiv0$ if $c>1$, and for every pair if $c=1$.
The endpoint can also be expressed pointwise with $G_2(z,w)=\Rea(zw)$ and $b_{2,c}=2a_{2,c}^2$.
\\
Corollary~\ref{cor:projection} completes the
determination of the sharp weighted
analytic--coanalytic projection norm for
all $1<p<\infty$.
\\
An additional feature of the full-range formula
is its invariance under conjugate exponents.
Indeed, if $p'=p/(p-1)$, then
\[
 a_{p,c}=a_{p',c},
 \qquad 1<p<\infty.
\]
This symmetry follows directly from the
explicit expression \eqref{eq:full-range};
the proof in the subquadratic range does not
rely on a scalar-to-vector duality identity.\\

\noindent\textbf{Funding.}
The authors gratefully acknowledge financial support from the
Ministry of Education, Science and Innovation of Montenegro through the
grants \emph{``Mathematical Analysis, Optimisation and Machine Learning''}
and \emph{``Complex-analytic and geometric techniques for non-Euclidean
machine learning: theory and applications.''}

\end{document}